# De Divino Errore

*'De Divina Proportione' was written by Luca Pacioli and illustrated by Leonardo da Vinci. It was one of the most widely read mathematical books. Unfortunately, a strongly emphasized statement in the book claims six summits of pyramids of the stellated icosidodecahedron lay in one plane. This is not so, and yet even extensively annotated editions of this book never noticed this error. Dutchmen Jos Janssens and Rinus Roelofs did so, 500 years later.*

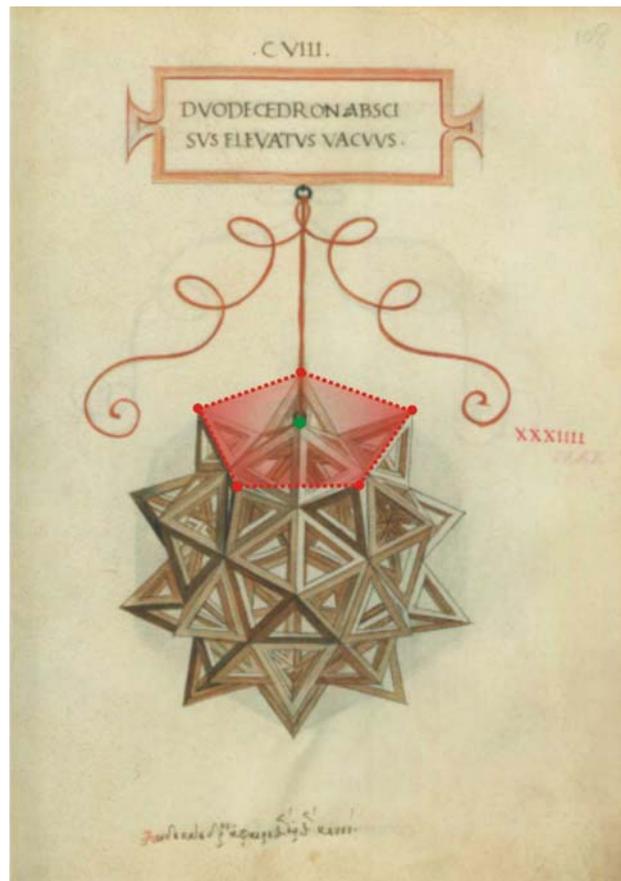

*Fig. 1: About this illustration of Leonardo da Vinci for the Milanese version of the 'De Divina Proportione', Pacioli erroneously wrote that the red and green dots lay in a plane.*

The book 'De Divina Proportione', or 'On the Divine Ratio', was written by the Franciscan Fra Luca Bartolomeo de Pacioli (1445-1517). His name is sometimes written Paciolo or Paccioli because Italian was not a uniform language in his days, when, moreover, Italy was not a country yet. Labeling Pacioli as a Tuscan, because of his birthplace of Borgo San Sepolcro, may be more correct, but he also studied in Venice and Rome, and spent much of his life in Perugia and Milan. In service of Duke and patron Ludovico Sforza, he would write his masterpiece, in 1497 (although it is more correct to say the work was written between 1496 and 1498, because it contains several parts). It was not his first opus, because in 1494 his 'Summa de arithmetic, geometrica, proportioni et proportionalita' had appeared; the 'Summa' and 'Divina' were not his only books, but surely the most famous ones.

For hundreds of years the books were among the most widely read mathematical bestsellers, their fame being only surpassed by the 'Elements' of Euclid. The 'Summa' was more popular than the 'Divina', though this no longer applies today because we now are especially impressed by the illustrations Leonardo da Vinci made for this book. They were the first truly insightful spatial representations of polyhedra, similar to contemporary 3D computer drawings. Perhaps, this collaboration between Pacioli and Leonardo was glorified in a well-known portrait from 1495 (see Fig. 2). It is usually attributed to Jacopo de' Barbari but this statement is questionable, as is the suggestion Leonardo was indeed the person painted in the background.

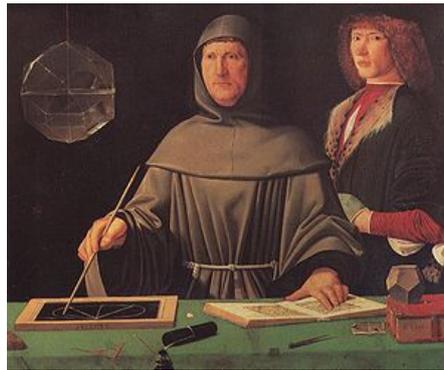

*Fig. 2: Portrait (1495) of Luca Pacioli and (?) Leonardo da Vinci.*

Pacioli's five hundred year old text seems tedious and unreadable to modern standards. Besides, the question arises whether 'The Divina Proportione' was indeed well read and verified, because chemist Jos Janssens (Leiderdorp, The Netherlands) recently happened to notice a curious statement in the book. Of some chapters about the polyhedra Pacioli studied, Janssens did a careful reading, inspired by recent discoveries of errors in Leonardo's geometric drawings (see [1], [2] and [3]). He found a questionable statement, printed in black and white, and thus is cannot be rejected as a matter of interpretation or as an inaccuracy. Janssens corresponded about this statement with the author and with his compatriot Rinus Roelofs (Hengelo), who found out Pacioli's statement was not correct. However, it is difficult to imagine that the error was never noticed, since there are so many publications about this work. However, they are mainly art history books, with a very limited focus on mathematics. Hence the current article about this error – though this paper can also be seen as a call to find out whether the age-old mistake is known or not.

**An incorrect statement**
Because of financial difficulties, 'The Divina Proportione' was only printed in 1509, in Venice, when Paganinus de Paganinus decided to take up the challenge of editing a mathematical book. Thus, it was more than ten years after the publication of the manuscript, of which two versions remain today, one in Geneva and one in Milan. In the latter two, the artwork by Leonardo da Vinci is hand-painted while in the printed version they were replaced by woodcuts. Their accuracy is a subject of debate ([4]) and this is good to know, in order to estimate precision and accuracy of evidence in those times. In any case, in the current paper the discussion is not about a drawing, but a statement that was clearly written down in precise words.

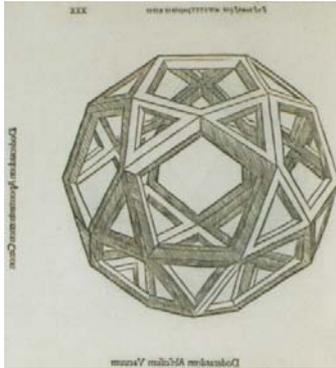

*Fig. 3: The truncated dodecahedron in the printed version of 'De Divina Proportione'.*

The statement is about the stellated icosidodecahedron, which is built on a truncated dodecahedron (see Fig. 3), as described in section 'LII' of 'The Divina Proportione' (the text was translated based on the French translation, see [4], and checked in the German version):

> *And the body that is created in this way, is composed by the flat truncated dodecahedron, inside, which shows itself to the mind only through the imagination, and by 32 pyramids, of which 12 are pentagonal, all of equal height, and of which the 20 others are triangular, all of equal height. The bases of the pyramids are the faces of the aforementioned dodecahedron and they mutually correspond, that is to say, the triangles to the triangular pyramids and the pentagons to the pentagonal pyramids. Projected onto a plane, this body will always rest on 6 tops of pyramids, one of them being a pentagonal pyramid, the other five triangular. When this body is seen in the air, it seems at first sight absurd that these vertices satisfy this property, but something like this, noble Duke, is of such great abstractness and deep science that I know that who understands me will not deny it. As for the dimensions of this body, they are obtained by the very subtle practice of primarily algebra or almucabala, which is known to few people, and is well demonstrated by us in our work, with methods allowing understanding them easily.*

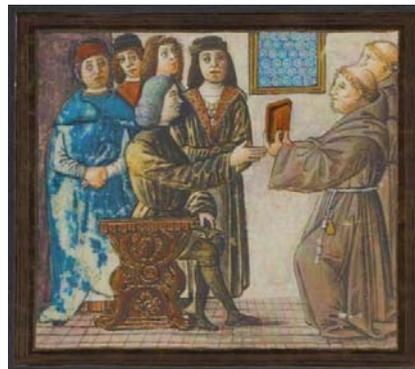

*Fig. 4: Fra Luca Pacioli donating his book to Duke Ludovico Sforza.*

Addressing the 'noble Duke ' in such a mathematical text was not unusual in those times. After all, Ludovico Sforza was the patron who had paid Pacioli (see Fig. 4). The self-flattering wording about the use of 'algebra and almucabala' seems outdated too, though the author right when he pretends checking geometric properties through calculation was little known in his time – after all, Descartes was just born. Now it happened Jos Janssens did what Pacioli said; he looked at the stellated icosidodecahedron. It seemed at first glance that the six points mentioned by Pacioli (the five red and one green on Fig. 1) indeed lie in one plane. Janssens Pacioli wondered if this had been verified on a concrete physical model, or whether an 'algebra and almucabala' calculation had actually been performed.

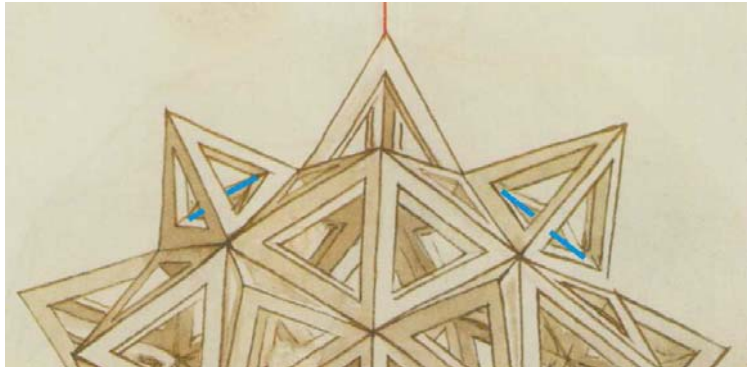
*Fig. 5: The same drawing as in Fig. 1, but from the Geneva version:
the blue lines illustrate an error in the drawing.*

Janssens immediately doubted about the geometric verification on a model, because the drawing of such a model showed a lack of precision which would not have occurred when the sketch would have been made from an existing model: some lines at the top of the figure are broken, while they should have been straight (see the blue lines on fig. 5). Because the trigonometric calculations seemed rather unattractive, Janssens asked Rinus Roelofs to check it out. Roelofs would not really calculate the statement either and do the algebraic work indirectly, using the modern 3D software Rhinosceros. And perhaps the lack of such software also explains why something was not observed during five hundred years, but now noticed by Janssens and Roelofs: the claim is not true. The green dot clearly comes out of the plane defined by the five red dots (see Fig. 6).

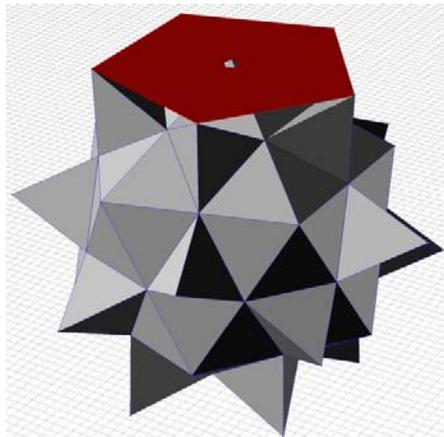
*Fig. 6: Roelofs ' negative answer to Janssens' question about Pacioli's theorem.*

Did Pacioli actually prove his 'theorem' about the six coplanar points, as he announced with so much commotion to Duke Sforza? A real model he probably did not have, otherwise he would have found out the polyhedron wobbles when standing on the top of a five-sided pyramid. Why no one ever noticed an error in such a popular book may be easier to explain: the calculations are quite complicated (though not impossible) and long-winding and perhaps no one took time and effort to actually execute them. Besides, Roelofs also let the calculations to a computer. Or maybe someone noticed the error, but it got lost in the numerous historic and artistic considerations that have glorified the work of Pacioli and Leonardo - rightfully so.

**References**

[1] Huylebrouck. Dirk (2011). '*Lost in Triangulation: Leonardo da Vinci's Mathematical Slip-Up*', http://www.scientificamerican.com/article.cfm?id=davinci-mathematical-slip-up.
[2] Huylebrouck. Dirk (2012). '*Lost in Enumeration: Leonardo da Vinci's Slip-Ups in Arithmetic and Mechanics'.* The Mathematical Intelligencer; DOI: 10.1007/s00283-012-9326-8.
[3] Huylebrouck. Dirk (2013). *'Lost in Edition: did Leonardo da Vinci slip up?'.* To appear in *Leonardo, 2014*.
[4] Luca Pacioli (1991). '*Divine Proportion'.* Librairie du Compagnonnage, 1980.